\RequirePackage{fix-cm}
\documentclass[smallextended]{svjour3}       % onecolumn (second format)
\smartqed  % flush right qed marks, e.g. at end of proof
\usepackage{graphicx}
\usepackage{bm}
\usepackage{fancyhdr}
\usepackage{graphicx}
\usepackage{subfigure}
\usepackage{amssymb}
\usepackage{amsmath}
\usepackage{amsfonts}
\usepackage{mathrsfs}
\usepackage{mathtools}

\usepackage{color}
\usepackage{setspace}
\usepackage{exscale}
\usepackage{hyperref}
\usepackage{relsize}
\usepackage{epstopdf}
\usepackage{float}
\usepackage{cite}

\usepackage{algorithm}
\usepackage{algpseudocode}
\usepackage{color}

\usepackage{hyperref}
\usepackage{setspace}
\usepackage{placeins}
\usepackage{booktabs,multirow} % for much better looking tables
\usepackage{array} % for better arrays (eg matrices) in maths
\usepackage{paralist} % very flexible & customisable lists (eg. enumerate/itemize, etc.)
\usepackage{verbatim} % adds environment for commenting out blocks of text & for better verbatim
\usepackage{subfigure} % make it possible to include more than one captioned figure/table in a single float
\allowdisplaybreaks[1]

\numberwithin{equation}{section}
\numberwithin{figure}{section}
\numberwithin{table}{section}

\newcommand\eref[1]{(\ref{#1})}

\newcommand*\xbar[1]{%
  \hbox{%
    \vbox{%
      \hrule height 0.5pt % The actual bar
      \kern0.4ex%         % Distance between bar and symbol
      \hbox{%
        \kern-0.05em%      % Shortening on the left side
        \ensuremath{#1}%
        \kern-0.00em%      % Shortening on the right side
      }%
    }%
  }%
}

\usepackage{xcolor}

\newcommand{\mF}{\bm{F}}

\newcommand{\mG}{\bm{G}}

\newcommand{\mU}{\bm{U}}

\newcommand{\mo}{\bm{0}}

\newcommand{\dx}{\Delta x}
\newcommand{\dy}{\Delta y}

\newcommand{\hf}{{\frac{1}{2}}}

\newcommand{\jph}{{j+\frac{1}{2}}}
\newcommand{\jmh}{{j-\frac{1}{2}}}
\newcommand{\kph}{{k+\frac{1}{2}}}
\newcommand{\kmh}{{k-\frac{1}{2}}}

\begin{document}

\title{Entropy-Based Local Characteristic Decomposition \thanks{The work of S. Chu and M. Herty was funded by the Deutsche
Forschungsgemeinschaft (DFG, German Research Foundation)--SPP 2410 Hyperbolic Balance Laws in Fluid Mechanics: Complexity, Scales,
Randomness (CoScaRa) within the Project HE5386/27-2 (Zuf\"allige kompressible Euler Gleichungen: Numerik und ihre Analysis, 525853336). The
work of A. Kurganov was supported in part by NSFC grant W2431004.}}

\titlerunning{Novel LCD}        % if too long for running head
\author{Shaoshuai Chu \and Michael Herty \and Alexander Kurganov}

\authorrunning{S. Chu, M. Herty \& A. Kurganov} % if too long for running head
\institute{S. Chu \at
Department of Mathematics, RWTH Aachen University, Aachen, 52056, Germany\\
\email{chu@igpm.rwth-aachen.de}%\\
%\emph{Present address:} of F. Author  %  if needed
\and
M. Herty \at
Department of Mathematics, RWTH Aachen University, 52056 Aachen, Germany; Department of Mathematics and Applied Mathematics, University of
Pretoria, Hatfield, 0028, South Africa\\
\email{herty@igpm.rwth-aachen.de}%\\
\and
A. Kurganov \at
Department of Mathematics and Shenzhen International Center for Mathematics, Southern University of Science and Technology, Shenzhen,
518055, China\\
\email{alexander@sustech.edu.cn}%\\
}

\date{Received: date / Accepted: date}

\maketitle

%\vspace{-1cm}
\begin{abstract}
Local characteristic decomposition (LCD) is widely used in high-order numerical methods for hyperbolic systems of conservation laws to
reduce spurious oscillations appearing in the computed solutions. The LCD implementation requires a representative average interface state,
typically obtained using arithmetic or Roe-type averages of the nearly grid values. We propose an entropy-based LCD (ELCD), in which the
states in the two cells adjacent to an interface are considered together with several nearby states in phase space. Among these candidates,
we select the state that locally minimizes the entropy and use it as an average interface state for linearizing the flux Jacobian. We
incorporate the ELCD into several second-order finite-volume and fifth-order finite-difference schemes. Numerical experiments for the
two-dimensional Euler equations of gas dynamics show that the schemes, which utilize the ELCD procedure generally resolve complex wave
structures more sharply than their counterparts, which use the arithmetic averages.


\end{abstract}

\section{Introduction}\label{intro}
This paper focuses on finite-volume (FV) and finite-difference (FD) numerical methods for hyperbolic systems of conservation laws, which, in
the two-dimensional (2-D) case, read as
\begin{equation}
\mU_t+\mF(\mU)_x+\mG(\mU)_y=\mo.
\label{1.2}
\end{equation}
Here, $t$ denotes the time, $x$ and $y$ are spatial variables, $\mU\in\mathbb R^m$ is the vector of conserved variables, and
$\mF,\mG:\mathbb R^m\to\mathbb R^m$ are the flux functions.

It is well-known that solutions of system \eref{1.2} may develop discontinuities in finite time even when the initial data are smooth. Such
solutions may contain complex wave structures involving shock waves, rarefaction waves, and contact discontinuities. The presence of these
nonsmooth structures makes the development of accurate and robust numerical methods particularly challenging. We refer the reader to the
monographs and review papers \cite{Hesthaven18,KLR20,Shu09,Shu20,Tor} and references therein for detailed discussions of various existing
numerical methods, including high-order ones.

A widely used framework for constructing FV and FD schemes for \eref{1.2} is based on a semi-discretization of \eref{1.2}. In this
framework, the flux derivatives are approximated using appropriate numerical fluxes, while the resulting system of time-dependent ODEs is
numerically integrated using an appropriate ODE solver. To achieve high-order spatial accuracy, one may proceed as follows. Given the
discrete data (either the cell averages or point values of the computed solution), one performs a piecewise polynomial reconstruction to
evaluate the one-sided point values of $\mU$ at the cell interfaces, where the numerical fluxes are to be computed. To make the
reconstruction non-oscillatory, the polynomial pieces must be constructed using a nonlinear limiter. We refer the reader to
\cite{Hesthaven18,KLR20,Shu09,Shu20,Tor}, where a variety of existing limiting approaches is discussed.

In the system case ($m>1$), high-order reconstructions are often performed in local characteristic variables to reduce spurious oscillations
that may appear when a componentwise reconstruction is used; see, e.g., \cite{Qiu02,Shu20}. To switch to the characteristic variables, we
perform the local characteristic decomposition (LCD) for the Jacobians linearized at the points, at which the numerical fluxes are to be
evaluated. A crucial component of the LCD is the construction of the representative average interface state. In numerical implementations,
this state is typically obtained using either the arithmetic or Roe-type averages of the neighboring solution values; see, e.g.,
\cite{Roe1981,Hesthaven18,Qiu02,Shu20}. However, it is not known what the best linearization option is. Though it is widely believed that
the choice of the average interface state almost does not affect the numerical solution, our studies suggest that the achieved resolution of
complex wave structures can be significantly improved by choosing an appropriate average interface state. 

In this paper, we propose a new entropy-based LCD (ELCD), for which the average interface state is selected as follows. We take the states
from the cells neighboring to the given cell interface as well as several states located nearby in the phase space. We then select the one
of those states, which locally minimizes the entropy, and make it to be the average interface state for the linearization of the Jacobian. 

We incorporate the proposed ELCD into second-order piecewise linear MinMod2 reconstruction \cite{lie03,Sweby84} and fifth-order WENO-Z
\cite{Borges08,wang18} interpolation performed in the corresponding local characteristic fields. We then test the new
reconstruction/interpolation in the second-order FV and fifth-order FD A-WENO \cite{JSZ,Liu17,CKX23} schemes, which are implemented using
four different numerical fluxes: the LLF \cite{KTcl,Rus61}, central-upwind (CU) from \cite{KNP}, low-dissipation CU (LDCU) \cite{CKX_24},
and HLLC \cite{TSS1994} ones. The resulting schemes are applied to the 2-D Euler equations of gas dynamics. The numerical results
demonstrate that for all of the studied schemes, the ELCD-based reconstruction/interpolation produces reliable results and often leads to
the improved resolution compared with the one achieved using the LCD based on the arithmetic average of the primitive variables.  

The rest of the paper is organized as follows. In \S\ref{sec2}, we overview the second-order FV and fifth-order FD A-WENO semi-discrete
schemes. In \S\ref{sec3}, we describe the proposed ELCD for the Euler equations of gas dynamics, which are then incorporated into the
piecewise polynomial reconstruction/interpolation and tested on a number of numerical examples in \S\ref{sec4}, where we illustrate the
potential advantages of the proposed ELCD.  Finally, we give some concluding remarks in \S\ref{sec5}.

\section{Semi-Discrete FV and FD Schemes: a Brief Overview}\label{sec2}
In this section, we introduce the second-order FV and fifth-order FD A-WENO schemes for the 2-D Euler equations of gas dynamics, which read
as \eref{1.2} with $\mU=(\rho,\rho u,\rho v,E)^\top$, $\mF=\big(\rho u,\rho u^2+p,\rho uv,u(E+p)\big)^\top$, and
$\mG=\big(\rho v,\rho uv,\rho v^2+p,v(E+p)\big)^\top$. Here, $\rho$, $u$, $v$, $p$, and $E$ denote the density, $x$- and $y$-directional
velocity, pressure, and total energy, respectively. The system is closed by the equation of state (EOS) for ideal gases:
\begin{equation}
p=(\gamma-1)\Big[E-\frac{\rho}{2}(u^2+v^2)\Big],
\label{2.1}
\end{equation}
where $\gamma$ is the ratio of specific heats. The (mathematical) entropy function for the Euler equations of gas dynamics is
\begin{equation}
{\cal S}=-\frac{\rho}{\gamma-1}\ln\Big(\frac{p}{\rho^\gamma}\Big).
\label{3.1a}
\end{equation}

\subsection{Second-Order Semi-Discrete FV Schemes}
Let the computational domain be partitioned into uniform Cartesian cells $I_{j,k}:=[x_\jmh,x_\jph]\times[y_\kmh,y_\kph]$, centered at
$(x_j,y_k)$ with $x_\jph-x_j=x_j-x_\jmh\equiv\dx/2$ and $y_\kph-y_k=y_k-y_\kmh\equiv\dy/2$. Assume that the computed cell averages,
$$
\xbar\mU_{j,k}:\approx\frac{1}{\dx\dy}\int\limits_{I_{j,k}}\mU(x,y,t)\,{\rm d}x{\rm d}y,
$$
are available at a certain $t\ge0$ (here and below, we omit the dependence of the indexed quantities on time for the sake of brevity). They
are evolved in time by numerically solving the following system of ODEs:
\begin{equation}
\frac{{\rm d}\xbar{\mU}_{j,k}}{{\rm d}t}=-\frac{\bm{{\cal F}}^{\rm FV}_{\jph,k}-\bm{{\cal F}}^{\rm FV}_{\jmh,k}}{\dx}
-\frac{\bm{{\cal G}}^{\rm FV}_{j,\kph}-\bm{{\cal G}}^{\rm FV}_{j,\kmh}}{\dy}.
\label{2.3}
\end{equation}
Here, $\bm{{\cal F}}^{\rm FV}_{\jph,k}=\bm{{\cal F}}^{\rm FV}\big(\mU^-_{\jph,k},\mU^+_{\jph,k}\big)$ and
$\bm{{\cal G}}^{\rm FV}_{j,\kph}=\bm{{\cal G}}^{\rm FV}\big(\mU^-_{j,\kph},\mU^+_{j,\kph}\big)$ denote the FV numerical fluxes. There are
many different numerical fluxes readily available in the literature. In \S\ref{sec4}, we will test four of them: the LLF, CU, LDCU, and HLLC
ones. These fluxes are computed using the one-sided point values $\mU^\pm_{\jph,k}$ and $\mU^\pm_{j,\kph}$, which are obtained with the help
of an appropriate piecewise linear reconstruction. Our particular choice used in the numerical examples reported in \S\ref{sec4} is the
MinMod2 reconstruction applied to the local characteristic variables as it was done, e.g., in \cite[Appendix E]{CKK2026}. 

\subsection{Fifth-Order Semi-Discrete A-WENO Schemes}
The FV second-order schemes can be extended to the fifth order of accuracy within the A-WENO framework, in which the point values
$\mU_{j,k}:\approx\mU(x_j,y_k,t)$ are evolved in time by numerically solving the following system of ODEs:
\begin{equation}
\frac{{\rm d}\mU_{j,k}}{{\rm d}t}=-\frac{\bm{{\cal H}}_{\jph,k}^x-\bm{{\cal H}}_{\jmh,k}^x}{\dx}
-\frac{\bm{{\cal H}}_{j,\kph}^y-\bm{{\cal H}}_{j,\kmh}^y}{\dy}.
\label{2.4}
\end{equation}
Here, $\bm{{\cal H}}_{\jph,k}^x$ and $\bm{{\cal H}}_{j,\kph}^y$ are the fifth-order A-WENO numerical fluxes 
\begin{equation*}
\begin{aligned}
&\bm{{\cal H}}_{\jph,k}^x=\bm{{\cal F}}^{\rm FV}_{\jph,k}-\frac{1}{24}(\dx)^2(\mF_{xx})_{\jph,k}+\frac{7}{5760}(\dx)^4(\mF_{xxxx})_{\jph,k},
\\
&\bm{{\cal H}}_{j,\kph}^y=\bm{{\cal G}}^{\rm FV}_{j,\kph}-\frac{1}{24}(\dy)^2(\mG_{yy})_{j,\kph}+\frac{7}{5760}(\dy)^4(\mG_{yyyy})_{j,\kph},
\end{aligned}
\end{equation*}
where $\bm{{\cal F}}^{\rm FV}_{\jph,k}$, $\bm{{\cal G}}^{\rm FV}_{j,\kph}$ are FV numerical fluxes and $(\mF_{xx})_{\jph,k}$,
$(\mF_{xxxx})_{\jph,k}$, $(\mG_{yy})_{j,\kph}$, $(\mG_{yyyy})_{j,\kph}$ are higher-order correction terms, computed using central
differences, which can be applied to already available FV numerical fluxes; see, e.g., \cite{CKX23}.

\iffalse 
In \cite{CKX23}, these terms were calculated using the already available FV numerical fluxes
by setting
\begin{equation*}
\begin{aligned}
&(\mF_{xx})_{\jph,k}=\frac{1}{12(\dx)^2}\Big[-\bm{{\cal F}}^{\rm FV}_{j-\frac{3}{2},k}+16\bm{{\cal F}}^{\rm FV}_{\jmh,k}-
30\bm{{\cal F}}^{\rm FV}_{\jph,k}\\
&\hspace{3.7cm}+16\bm{{\cal F}}^{\rm FV}_{j+\frac{3}{2},k}-\bm{{\cal F}}^{\rm FV}_{j+\frac{5}{2},k}\Big],\\
&(\mF_{xxxx})_{\jph,k}=\frac{1}{(\dx)^4}\Big[\bm{{\cal F}}^{\rm FV}_{j-\frac{3}{2},k}-4\bm{{\cal F}}^{\rm FV}_{\jmh,k}+
6\bm{{\cal F}}^{\rm FV}_{\jph,k}-4\bm{{\cal F}}^{\rm FV}_{j+\frac{3}{2},k}+\bm{{\cal F}}^{\rm FV}_{j+\frac{5}{2},k}\Big],\\[0.5ex]
&(\mG_{yy})_{j,\kph}=\frac{1}{12(\dy)^2}\Big[-\bm{{\cal G}}^{\rm FV}_{j,k-\frac{3}{2}}+16\bm{{\cal G}}^{\rm FV}_{j,\kmh}-
30\bm{{\cal G}}^{\rm FV}_{j,\kph}\\
&\hspace{3.7cm}+16\bm{{\cal G}}^{\rm FV}_{j,k+\frac{3}{2}}-\bm{{\cal G}}^{\rm FV}_{j,k+\frac{5}{2}}\Big],\\
&(\mG_{yyyy})_{j,\kph}=\frac{1}{(\dy)^4}\Big[\bm{{\cal G}}^{\rm FV}_{j,k-\frac{3}{2}}-4\bm{{\cal G}}^{\rm FV}_{j,\kmh}+
6\bm{{\cal G}}^{\rm FV}_{j,\kph}-4\bm{{\cal G}}^{\rm FV}_{j,k+\frac{3}{2}}+\bm{{\cal G}}^{\rm FV}_{j,k+\frac{5}{2}}\Big].
\end{aligned}
\end{equation*}
\fi 

Finally, the one-sided point values $\mU^\pm_{\jph,k}$ and $\mU^\pm_{j,\kph}$ are obtained using an appropriate fifth-order interpolation.
In the numerical experiments, reported in \S\ref{sec4}, we have used the WENO-Z interpolation \cite{Borges08,wang18} applied to the
characteristic variables as it was done, e.g., in \cite[Appendix C]{CKK2026}.

\section{Entropy-Based Average Interface States}\label{sec3}
In this section, we describe how to select the average interface states and to use them in the ELCD reconstruction/interpolation procedure.
Here, we only show the details on the reconstruction/interpolation of $\mU^\pm_{\jph,k}$, as $\mU^\pm_{j,\kph}$ can be obtained in a similar
manner. 

We first introduce the $x$-directional Jacobian $A(\mU):=\frac{\partial\mF}{\partial\mU}$ and its linearization
$\widehat A_{\jph,k}=A\big(\widehat\mU_{\jph,k}\big)$, where $\widehat\mU_{\jph,k}$ is the entropy-based average interface state. We then
denote by $R_{\jph,k}=R\big(\widehat\mU_{\jph,k}\big)$ and $R^{-1}_{\jph,k}=R^{-1}\big(\widehat\mU_{\jph,k}\big)$ the corresponding right-
and left-eigenvector matrices, respectively, so that $R_{\jph,k}^{-1}\widehat A_{\jph,k}R_{\jph,k}$ is a diagonal matrix. Details on
$R_{\jph,k}$ and $R^{-1}_{\jph,k}$ for the 2-D Euler equations of gas dynamics can be found in \cite[Appendix C]{CCHKL_22}. 

We stress that both $R_{\jph,k}$ and $R^{-1}_{\jph,k}$ depend on the average interface state $\widehat\mU_{\jph,k}$. In
\cite{CCHKL_22,CKK2026} and many other works, these states were computed by arithmetic averages of the primitive variables from cells
$I_{j,k}$ and $I_{j+1,k}$, namely,
\begin{equation}
\begin{aligned}
&\widehat\rho_{\jph,k}=\frac{\,\xbar\rho_{j,k}+\,\xbar\rho_{j+1,k}}{2},\quad&\widehat u_{\jph,k}=\frac{u_{j,k}+u_{j+1,k}}{2},\\
&\widehat v_{\jph,k}=\frac{v_{j,k}+v_{j+1,k}}{2},\quad&\widehat p_{\jph,k}=\frac{p_{j,k}+p_{j+1,k}}{2},
\end{aligned}
\label{3.3a}
\end{equation}
where $u_{j,k}={(\xbar{\rho u})_{j,k}}/\,\xbar\rho_{j,k}$, $v_{j,k}={(\xbar{\rho v})_{j,k}}/\,\xbar\rho_{j,k}$, and $p_{j,k}$ is obtained by
the EOS \eref{2.1}:
\begin{equation}
p_{j,k}=(\gamma-1)\big[\,\xbar E_{j,k}-\hf\,\xbar\rho_{j,k}\big(u_{j,k}^2+v_{j,k}^2\big)\big].
\label{3.2f}
\end{equation}

This choice, however, is not optimal. We propose to use different, entropy-based average interface states, which are obtained by locally
minimizing the entropy within the rectangle in the $(\rho,p)$-plane with the following four vertices:
\begin{equation}
(\,\xbar\rho_{j,k},p_{j,k}),\quad(\,\xbar\rho_{j+1,k},p_{j,k}),\quad(\,\xbar\rho_{j,k},p_{j+1,k}),\quad{\rm and}\quad
(\,\xbar\rho_{j+1,k},p_{j+1,k}).
\label{3.2}
\end{equation}
In fact, since the entropy ${\cal S}$ defined in \eref{3.1a} is a convex function of $\rho$ and $p$, we only need to take the minimum
between the values of ${\cal S}$ at the four corner points in \eref{3.2}. The interface states are then determined by the following
algorithm.
\begin{algorithm}[ht!]
\caption{Entropy-based average interface state in the $x$-direction.}\label{alg2}
\begin{algorithmic}[1]
\State {\bf Set} ${\cal S}_1:={\cal S}(\,\xbar\rho_{j,k},p_{j,k})$, ${\cal S}_2:={\cal S}(\,\xbar\rho_{j+1,k},p_{j+1,k})$,

${\cal S}_3:={\cal S}(\,\xbar\rho_{j,k},p_{j+1,k})$, ${\cal S}_4:={\cal S}(\,\xbar\rho_{j+1,k}, p_{j,k})$
\vskip3pt
\If{${\cal S}_1=\min\{{\cal S}_1,{\cal S}_2,{\cal S}_3,{\cal S}_4\}$}
\State $\widehat\rho_{\jph,k}:=\,\xbar\rho_{j,k}$, $\widehat p_{\jph,k}:=p_{j,k}$, $\widehat u_{\jph,k}:=u_{j,k}$,
$\widehat v_{\jph,k}:=v_{j,k}$
\ElsIf{${\cal S}_2=\min\{{\cal S}_1,{\cal S}_2,{\cal S}_3,{\cal S}_4\}$}
\State $\widehat\rho_{\jph,k}:=\,\xbar\rho_{j+1,k}$, $\widehat p_{\jph,k}:=p_{j+1,k}$, $\widehat u_{\jph,k}:=u_{j+1,k}$,
$\widehat v_{\jph,k}:=v_{j+1,k}$
\ElsIf{${\cal S}_3=\min\{{\cal S}_1,{\cal S}_2,{\cal S}_3,{\cal S}_4\}$}
\State $\widehat\rho_{\jph,k}:=\,\xbar\rho_{j,k}$, $\widehat p_{\jph,k}:=p_{j+1,k}$,

$\hspace*{-0.15cm}\widehat u_{\jph,k}:=\hf(u_{j,k}+u_{j+1,k})$, $\widehat v_{\jph,k}:=\hf(v_{j,k}+v_{j+1,k})$
\Else
\State $\widehat\rho_{\jph,k}:=\,\xbar\rho_{j+1,k}$, $\widehat p_{\jph,k}:=p_{j,k}$,

$\hspace*{-0.15cm}\widehat u_{\jph,k}:=\hf(u_{j,k}+u_{j+1,k})$, $\widehat v_{\jph,k}:=\hf(v_{j,k}+v_{j+1,k})$
\EndIf
\end{algorithmic}
\end{algorithm}

\noindent
The selected primitive variables define the entropy-based conservative average interface state
\begin{equation*}
\widehat\mU_{\jph,k}=\left(\widehat\rho_{\jph,k},\widehat\rho_{\jph,k}\widehat u_{\jph,k},\widehat\rho_{\jph,k}\widehat v_{\jph,k},
\widehat E_{\jph,k}\right)^\top,
\end{equation*}
where
$$
\widehat E_{\jph,k}=\frac{\widehat p_{\jph,k}}{\gamma-1}+\frac{\widehat\rho_{\jph,k}}{2}
\left(\widehat u_{\jph,k}^{\,2}+\widehat v_{\jph,k}^{\,2}\right).
$$
\begin{remark}
For the fifth-order A-WENO schemes, the bars in \eref{3.3a}--\eref{3.2}, and Algorithm 1 should be removed since the point values
$\rho_{j,k}$, $(\rho u)_{j,k}$, $(\rho v)_{j,k}$, and $E_{j,k}$ are to be used there instead of the cell averages $\,\xbar\rho_{j,k}$,
$(\xbar{\rho u})_{j,k}$, $(\xbar{\rho v})_{j,k}$, and $\,\xbar E_{j,k}$. 
\end{remark}
\begin{remark}
We stress that Algorithm 1 is not completely symmetric with respect to the four candidates. Since the entropy $\cal S$ depends only on
$\rho$ and $p$, the minimization determines only $\widehat\rho_{\jph,k}$ and $\widehat p_{\jph,k}$. If ${\cal S}_1$ or ${\cal S}_2$ is the
minimum, the selected density and pressure are taken from the same cell, and we therefore use the velocities from that cell. If ${\cal S}_3$
or ${\cal S}_4$ is the minimum, the selected density and pressure are taken from different cells, and hence there are no corresponding
velocity values and we use the arithmetic averages of each of the velocity components.
\end{remark}

\section{Numerical Examples}\label{sec4}
In this section, we compare the performance of the second-order FV and fifth-order A-WENO schemes designed with the help of the ELCD and the
LCD based on \eref{3.3a}. We test four different numerical fluxes (the LLF, CU, LDCU, and HLLC ones) and denote the corresponding second-
and fifth-order schemes by LLF2, CU2, LDCU2, HLLC2 and LLF5, CU5, LDCU5, HLLC5, respectively. When the ELCD is used, we denote the resulting
schemes by LLF2-ELCD, CU2-ELCD, etc. 

We numerically integrate the semi-discrete systems \eref{2.3} and \eref{2.4} by the three-stage third-order strong stability preserving
Runge-Kutta (SSP RK3) method (see, e.g., \cite{Gottlieb11,Gottlieb12}) and use the CFL number $0.4$. In all of the numerical examples, we
take $\gamma=1.4$.

\subsubsection*{Example 1---2-D Riemann Problem (Configuration 12)}
In the first example, we consider Configuration 12 of the 2-D Riemann problems taken from \cite{Kurganov02}. The initial conditions,
\begin{equation*}
(\rho,u,v,p)\Big|_{(x,y,0)}=\begin{cases}(0.5313,0,0,0.4),&x>0.5,~y>0.5,\\(1,0.7276,0,1),&x<0.5,~y>0.5,\\(0.8,0,0,1),&x<0.5,~y<0.5,\\
(1,0,0.7276,1),&x>0.5,~y<0.5,\end{cases}
\end{equation*}
are prescribed in the computational domain $[0,0.6]\times[0,0.6]$ subject to the free boundary conditions.

\begin{figure}[ht!]
\centerline{\includegraphics[trim=0cm 0cm 0cm 0cm, clip, width=0.60\linewidth]{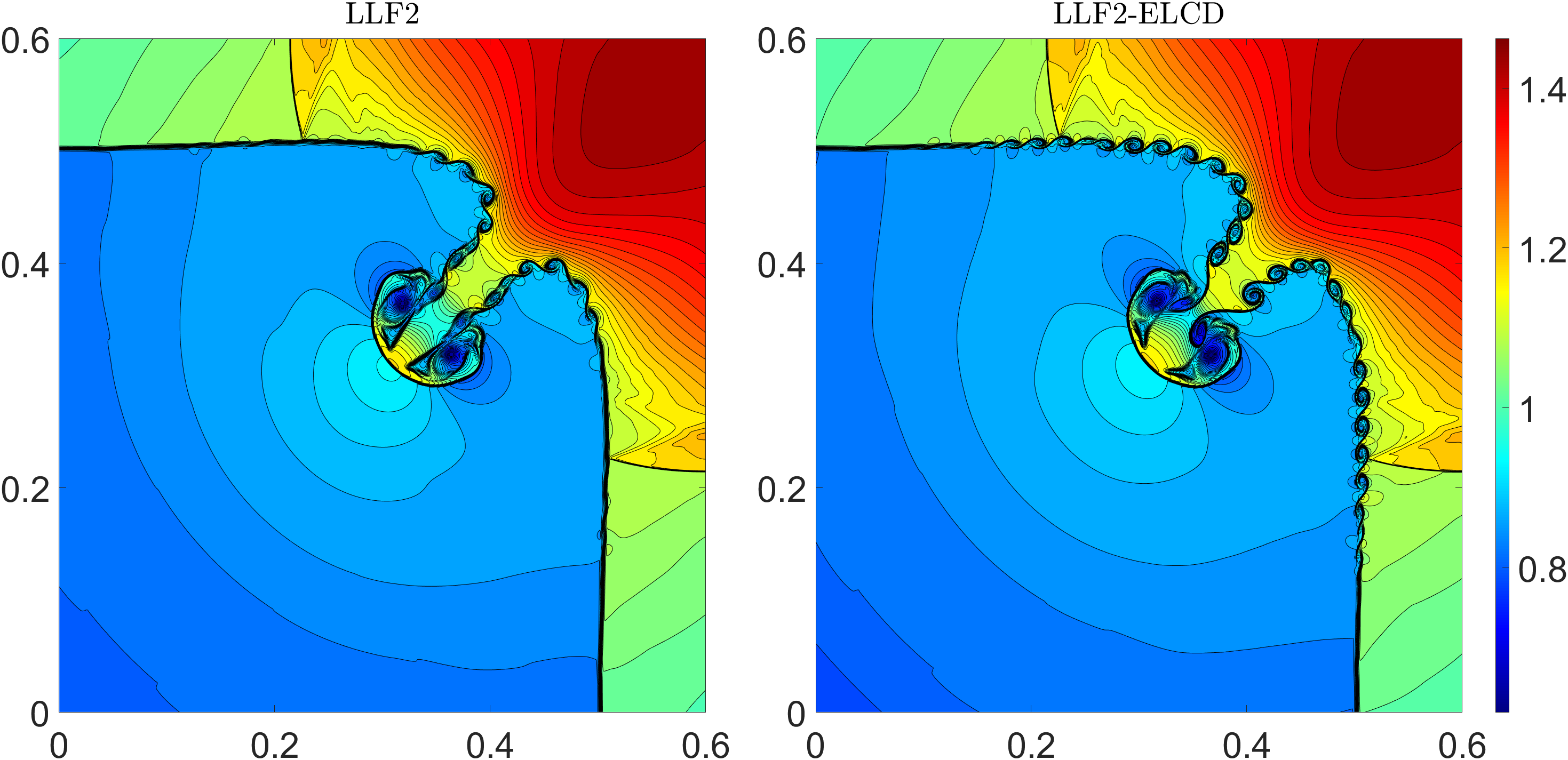}}
\vskip4pt
\centerline{\includegraphics[trim=0cm 0cm 0cm 0cm, clip, width=0.60\linewidth]{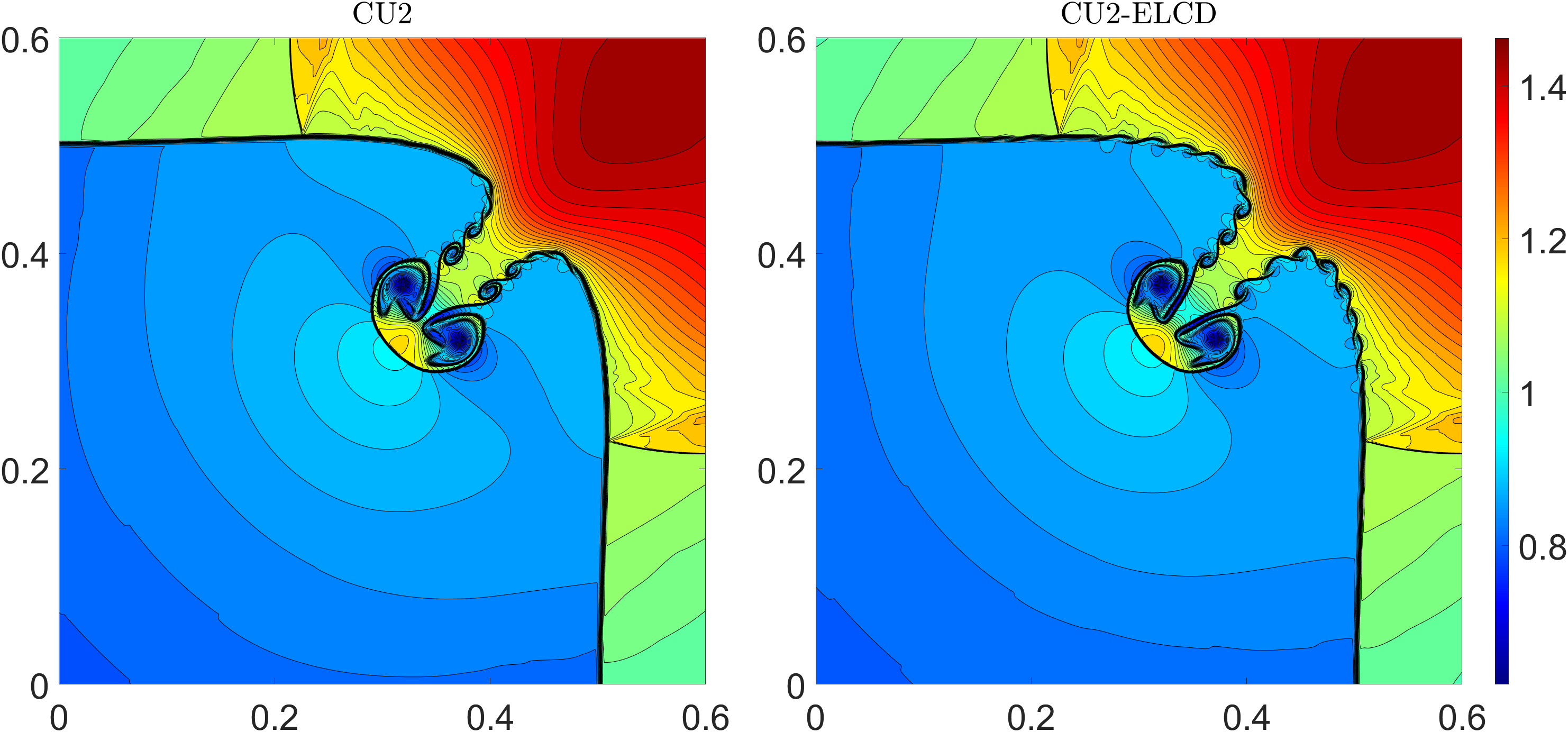}}
\vskip4pt
\centerline{\includegraphics[trim=0cm 0cm 0cm 0cm, clip, width=0.60\linewidth]{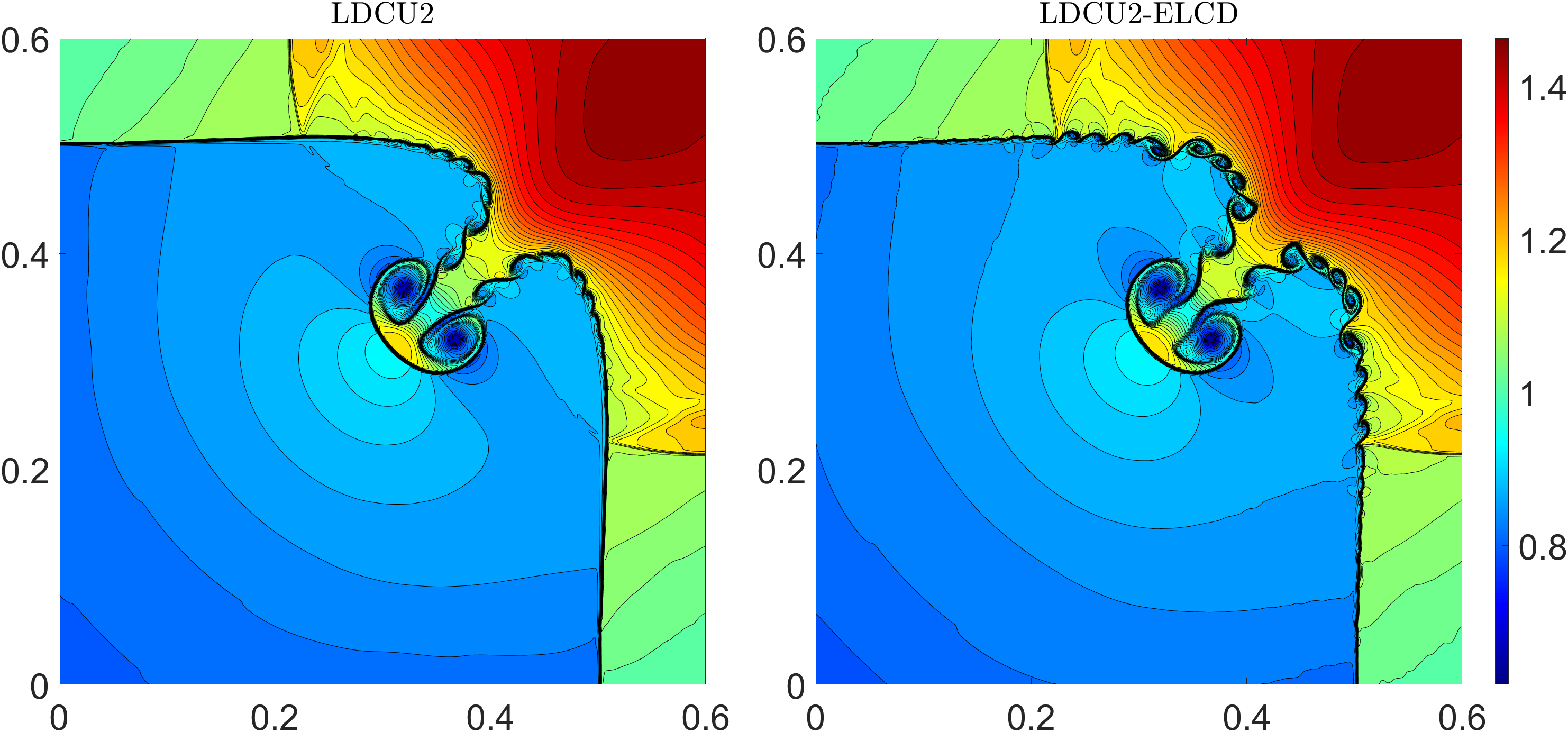}}
\vskip4pt
\centerline{\includegraphics[trim=0cm 0cm 0cm 0cm, clip, width=0.60\linewidth]{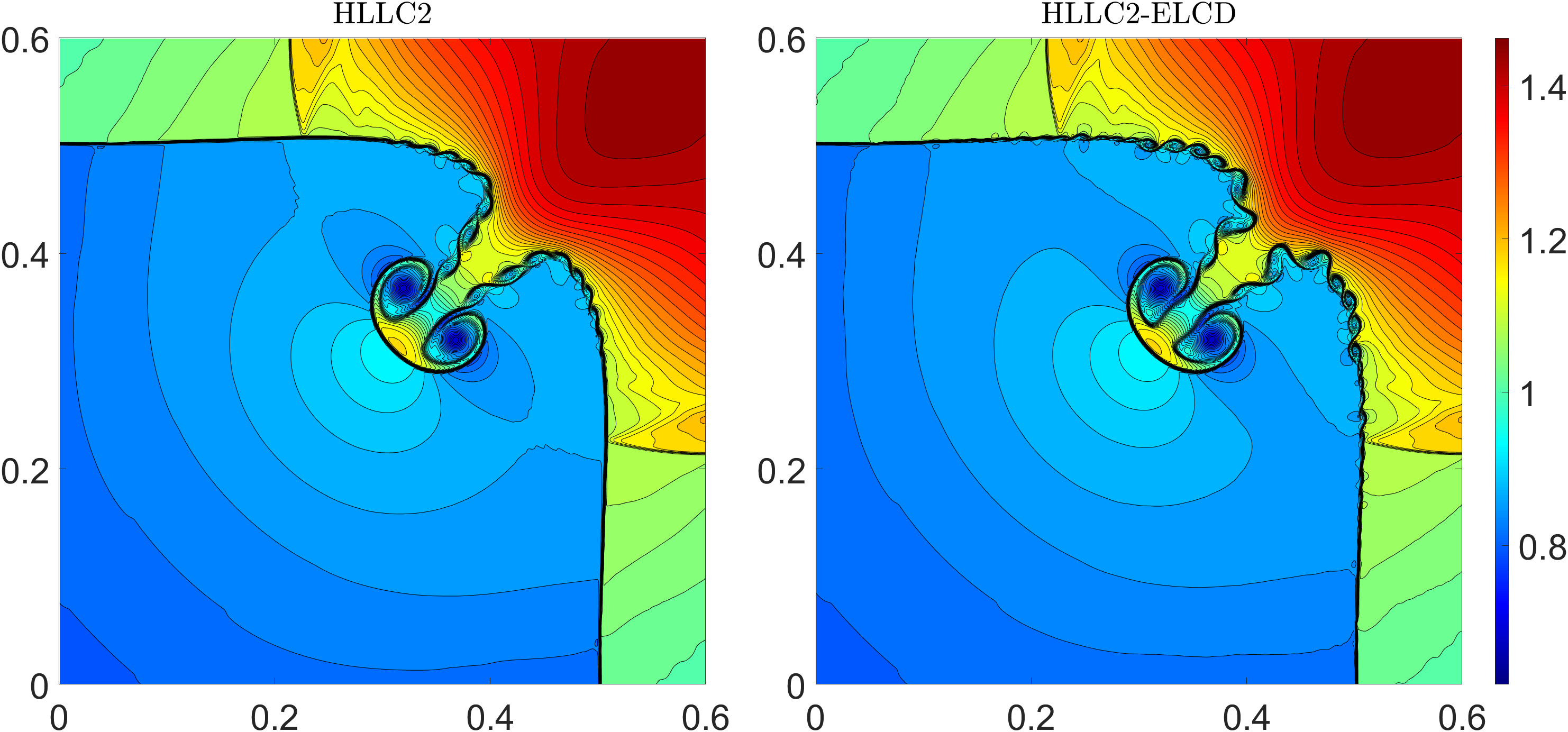}}
\caption{\sf Example 1 (second-order FV schemes): Density $\rho$ computed using the LLF (top row), CU (second row), LDCU (third row), and
HLLC (bottom row) numerical fluxes on $1800\times1800$, $1500\times1500$, $600\times600$, and $600\times600$ uniform meshes, respectively.
The results obtained by the schemes that use the ELCD procedure are plotted in the right column.\label{fig2}}
\end{figure}

\begin{figure}[ht!]
\centerline{\includegraphics[trim=0cm 0cm 0cm 0cm, clip, width=0.60\linewidth]{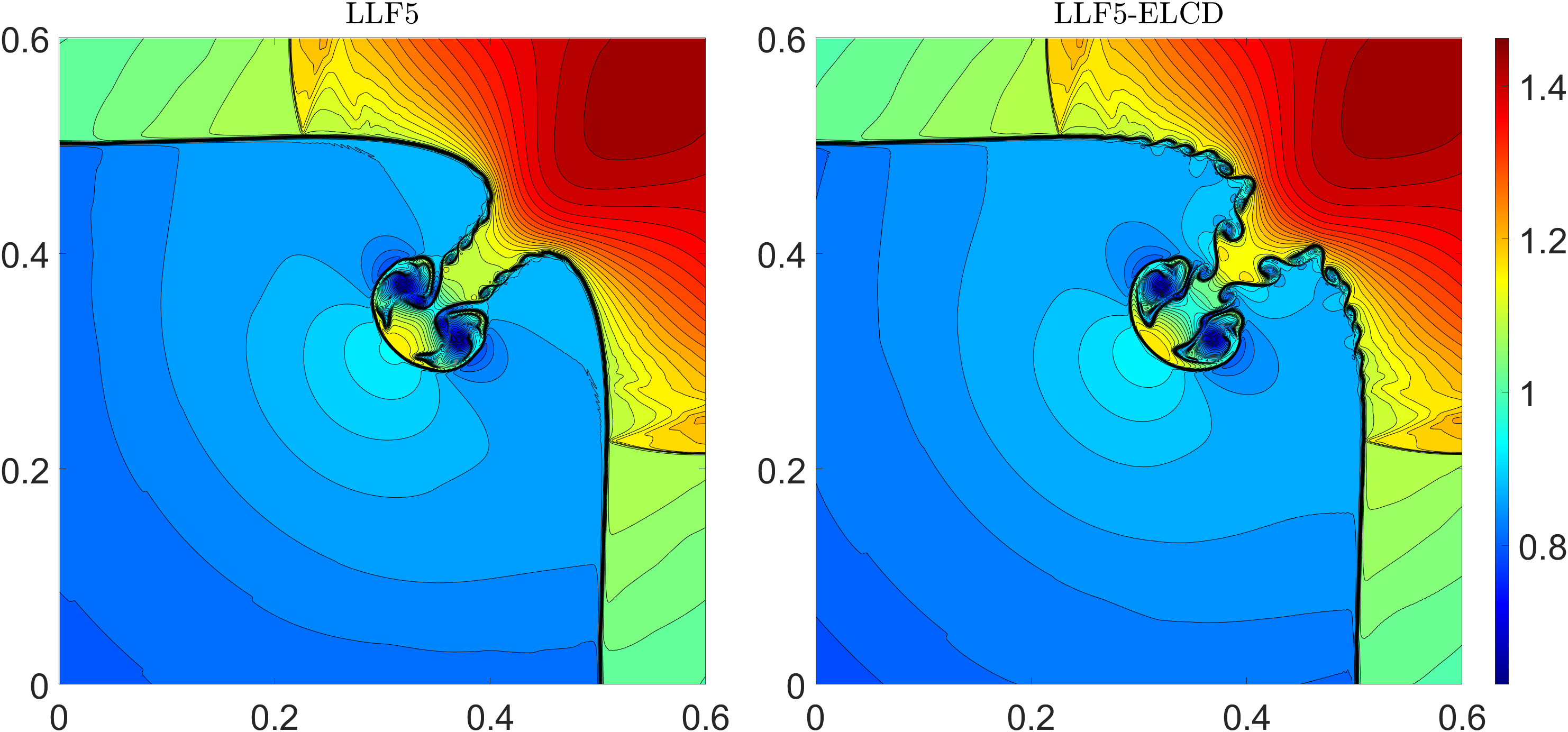}}
\vskip4pt
\centerline{\includegraphics[trim=0cm 0cm 0cm 0cm, clip, width=0.60\linewidth]{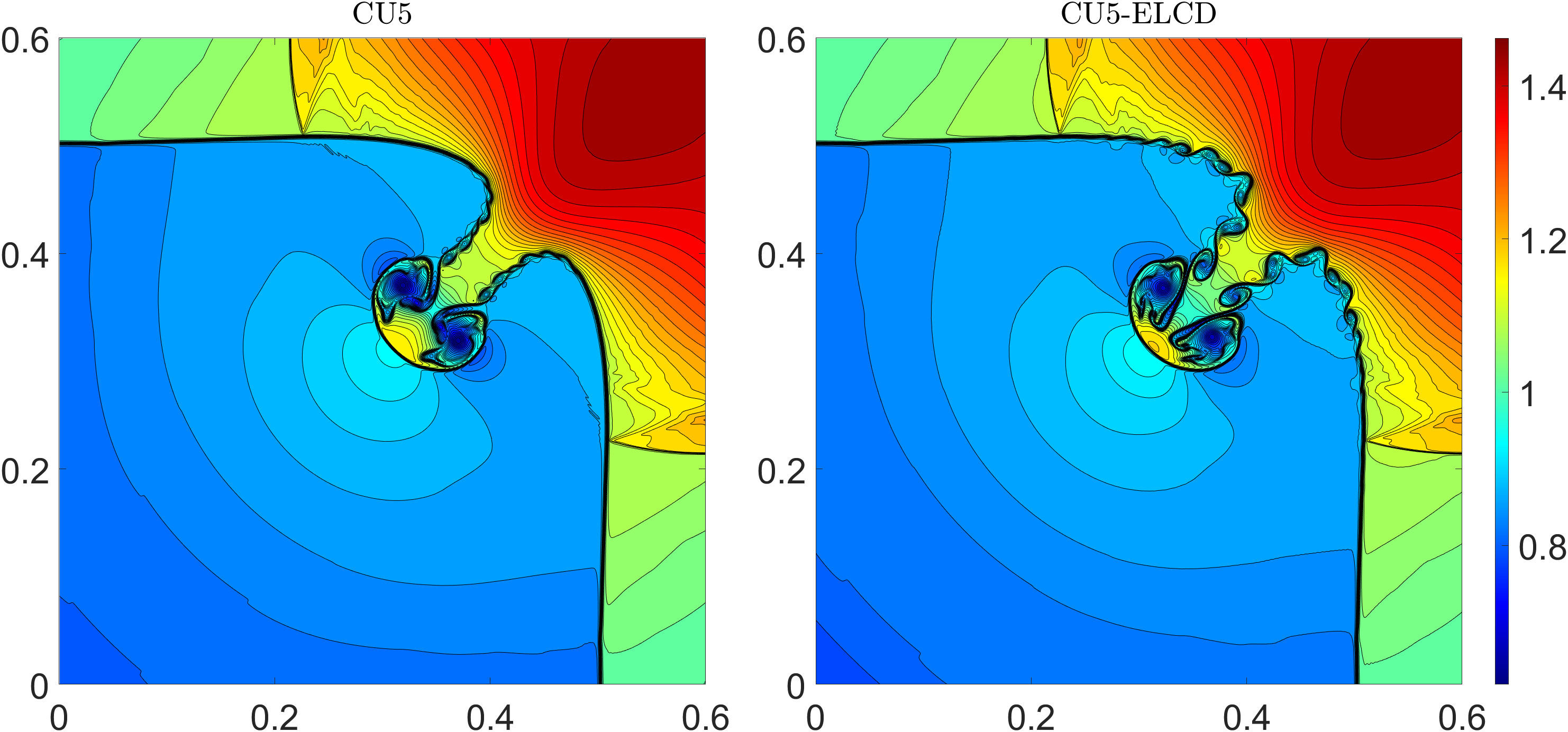}}
\vskip4pt
\centerline{\includegraphics[trim=0cm 0cm 0cm 0cm, clip, width=0.60\linewidth]{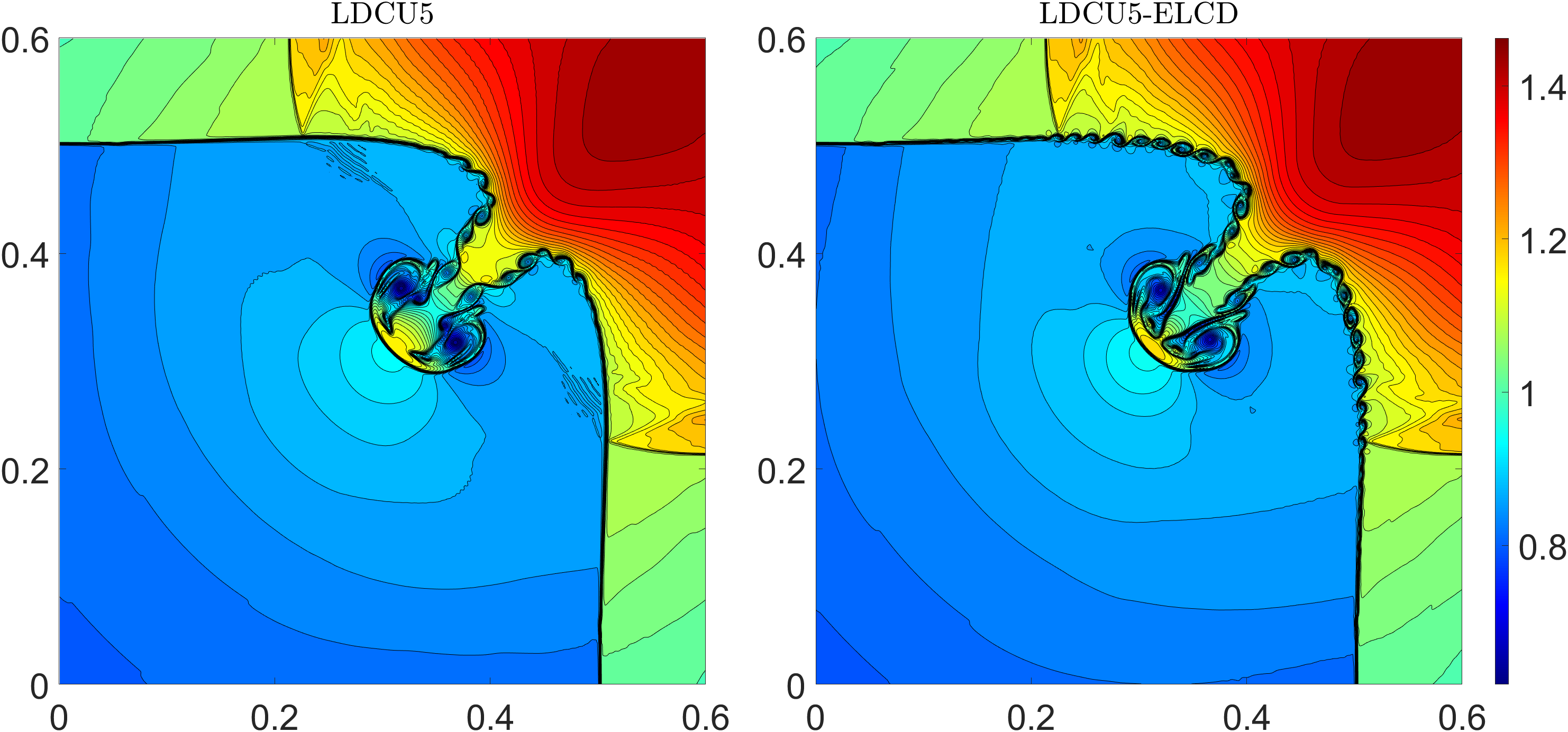}}
\vskip4pt
\centerline{\includegraphics[trim=0cm 0cm 0cm 0cm, clip, width=0.60\linewidth]{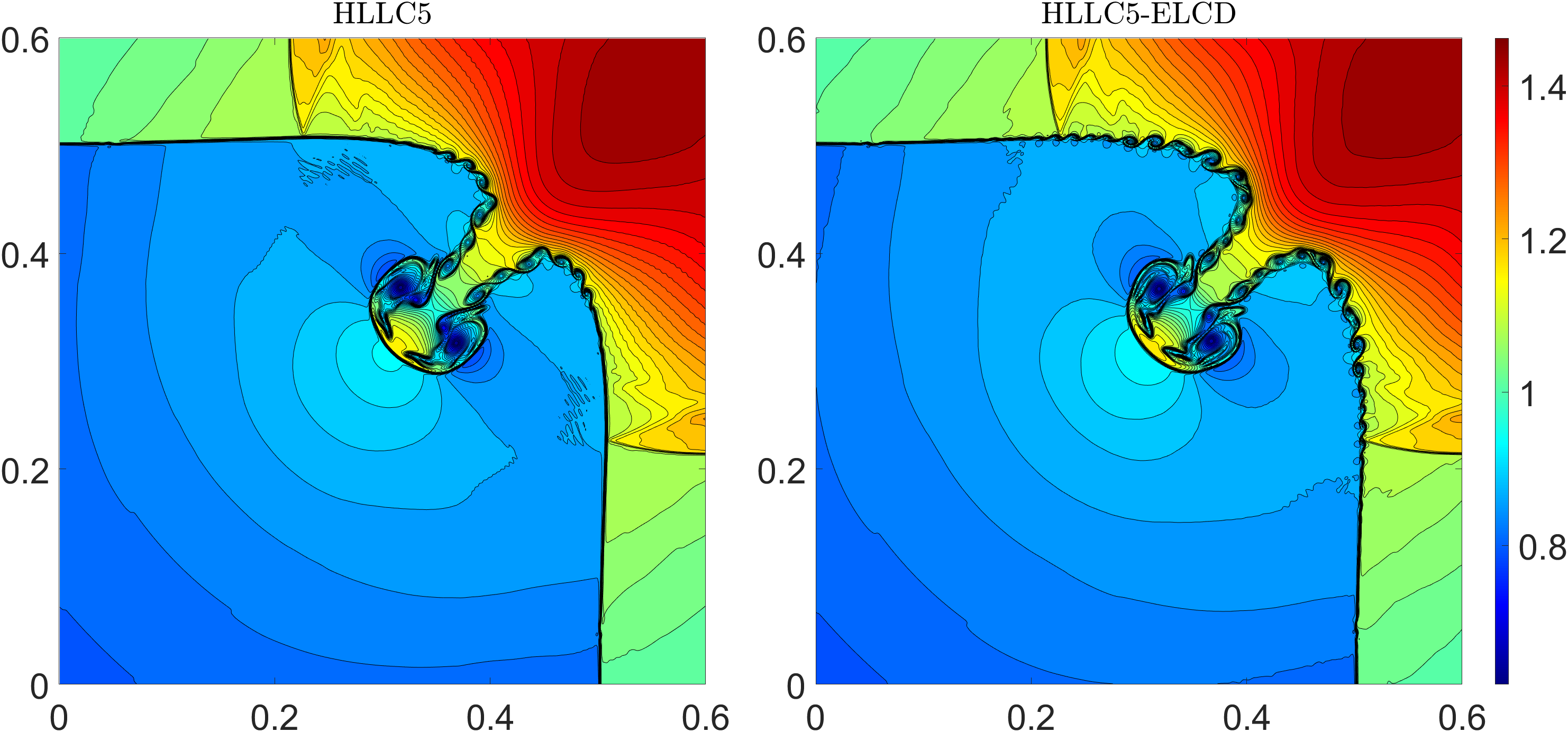}}
\caption{\sf Example 1 (fifth-order FD schemes): Density $\rho$ computed using the LLF (top row), CU (second row), LDCU (third row), and
HLLC (bottom row) numerical fluxes on $800\times800$, $800\times800$, $600\times600$, and $600\times600$ uniform meshes, respectively. The
results obtained by the schemes that use the ELCD procedure are plotted in the right column.\label{fig2a}}
\end{figure}

We compute the numerical solution until the final time $t=1$ by the studied second- and fifth-order schemes on different uniform meshes,
which are specified in the caption of Figs. \ref{fig2} and \ref{fig2a}, where the obtained results are presented. As one can see, the
schemes, which use the ELCD procedure, resolve finer structures including the vortices arising along the unstable contact surfaces compared
to their counterparts, which use the arithmetic averages \eref{3.3a}. 

\subsubsection*{Example 2--2-D Riemann Problem (Configuration 6)}
In this example, we consider Configuration 6 of the 2-D Riemann problems taken from \cite{Kurganov02}. The initial conditions,
\begin{equation*}
(\rho,u,v,p)(x,y,0)=\begin{cases}(1,0.75,-0.5,1),&x>0.5,~y>0.5,\\(2,0.75,0.5,1),&x<0.5,~y>0.5,\\(1,-0.75,0.5,1),&x<0.5,~y<0.5,\\
(3,-0.75,-0.5,1),&x>0.5,~y<0.5,\end{cases}
\end{equation*}
are prescribed in $[0,1]\times[0,1]$ subject to the free boundary conditions.

\begin{figure}[ht!]
\centerline{\includegraphics[trim=0cm 0cm 0cm 0cm, clip, width=0.60\linewidth]{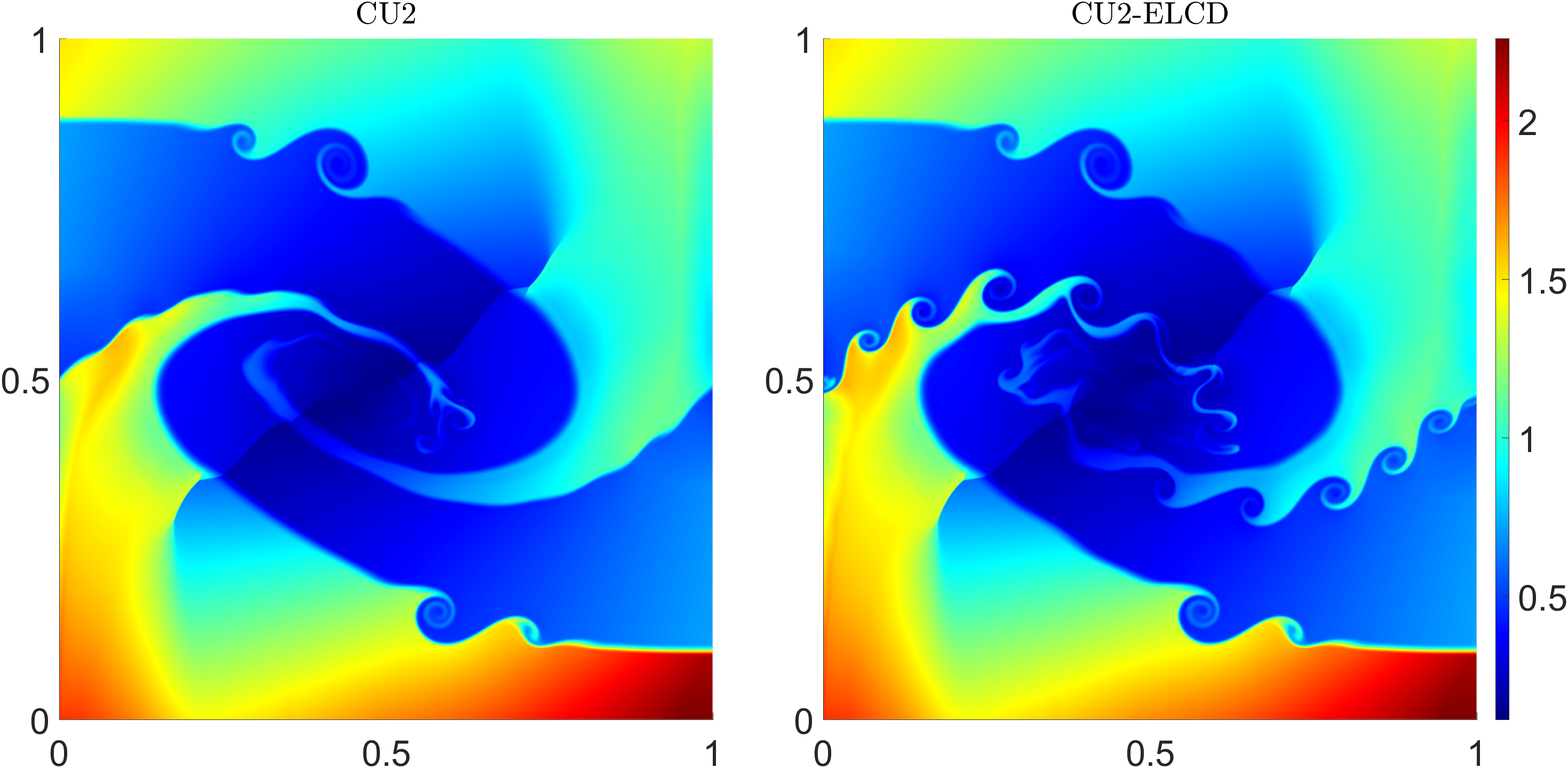}}
\vskip4pt 
\centerline{\includegraphics[trim=0cm 0cm 0cm 0cm, clip, width=0.60\linewidth]{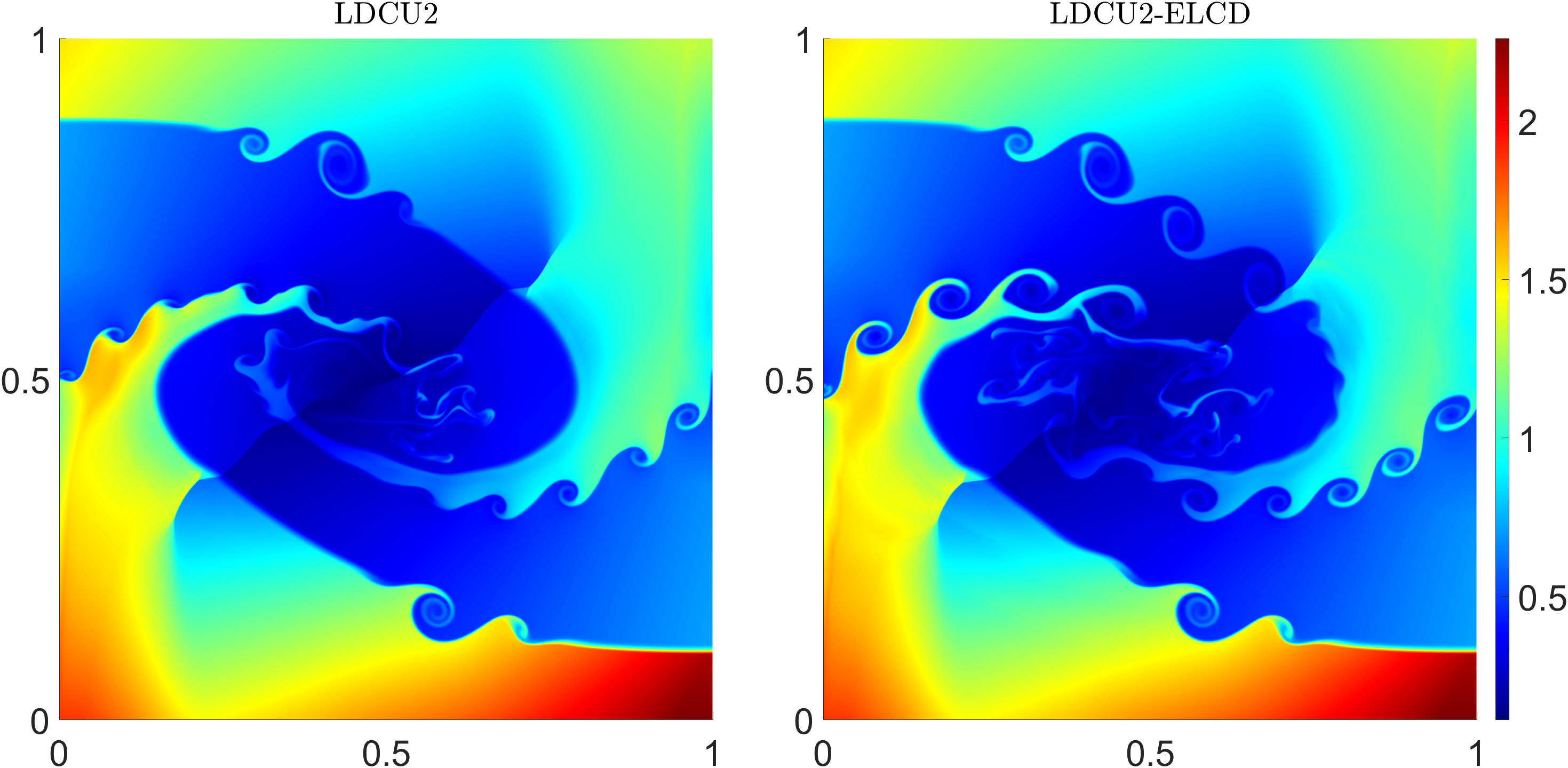}}
\vskip4pt 
\centerline{\includegraphics[trim=0cm 0cm 0cm 0cm, clip, width=0.60\linewidth]{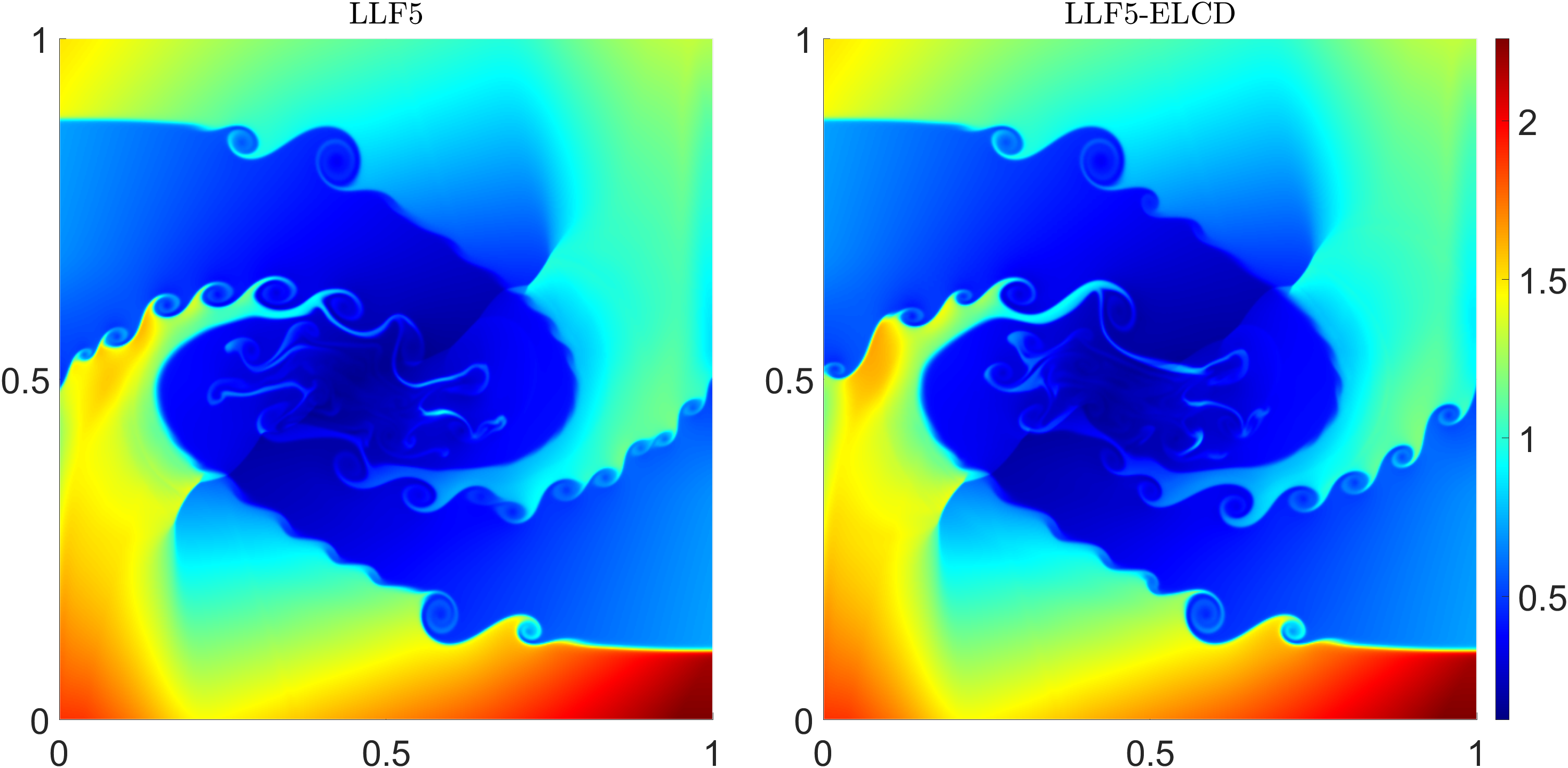}}
\caption{\sf Example 2: Density $\rho$ computed using the CU2 (top row), LDCU2 (middle row), and LLF5 (bottom row) schemes on
$1000\times1000$, $1000\times1000$, and $600\times600$ uniform meshes, respectively. The results obtained by the schemes that use the ELCD
procedure are plotted in the right column.\label{fig3}}
\end{figure}

We compute the numerical solution until the final time $t=1$ using the CU2, LDCU2, and LLF5 schemes on different uniform meshes specified in
the caption of Fig. \ref{fig3}. As one can see, the schemes, which use the ELCD procedure produce finer structures and hence achieve higher
resolution compared with their counterparts, which use the arithmetic averages \eref{3.3a}. 

\subsubsection*{Example 3---2-D Riemann Problem (Configuration 3)}
In this example, we consider Configuration 3 of the 2-D Riemann problems from \cite{Kurganov02}. The initial conditions,
\begin{equation*}
(\rho,u,v,p)(x,y,0)=\begin{cases}(1.5,0,0,1.5),&x>1,~y>1,\\(0.5323,1.206,0,0.3),&x<1,~y>1,\\(0.138,1.206,1.206,0.029),&x<1,~y<1,\\
(0.5323,0,1.206,0.3),&x>1,~y<1,\end{cases}
\end{equation*}
are prescribed in $[0,1.2]\times[0,1.2]$ subject to the free boundary conditions.

\begin{figure}[ht!]
\centerline{\includegraphics[trim=0cm 0cm 0cm 0cm, clip, width=0.60\linewidth]{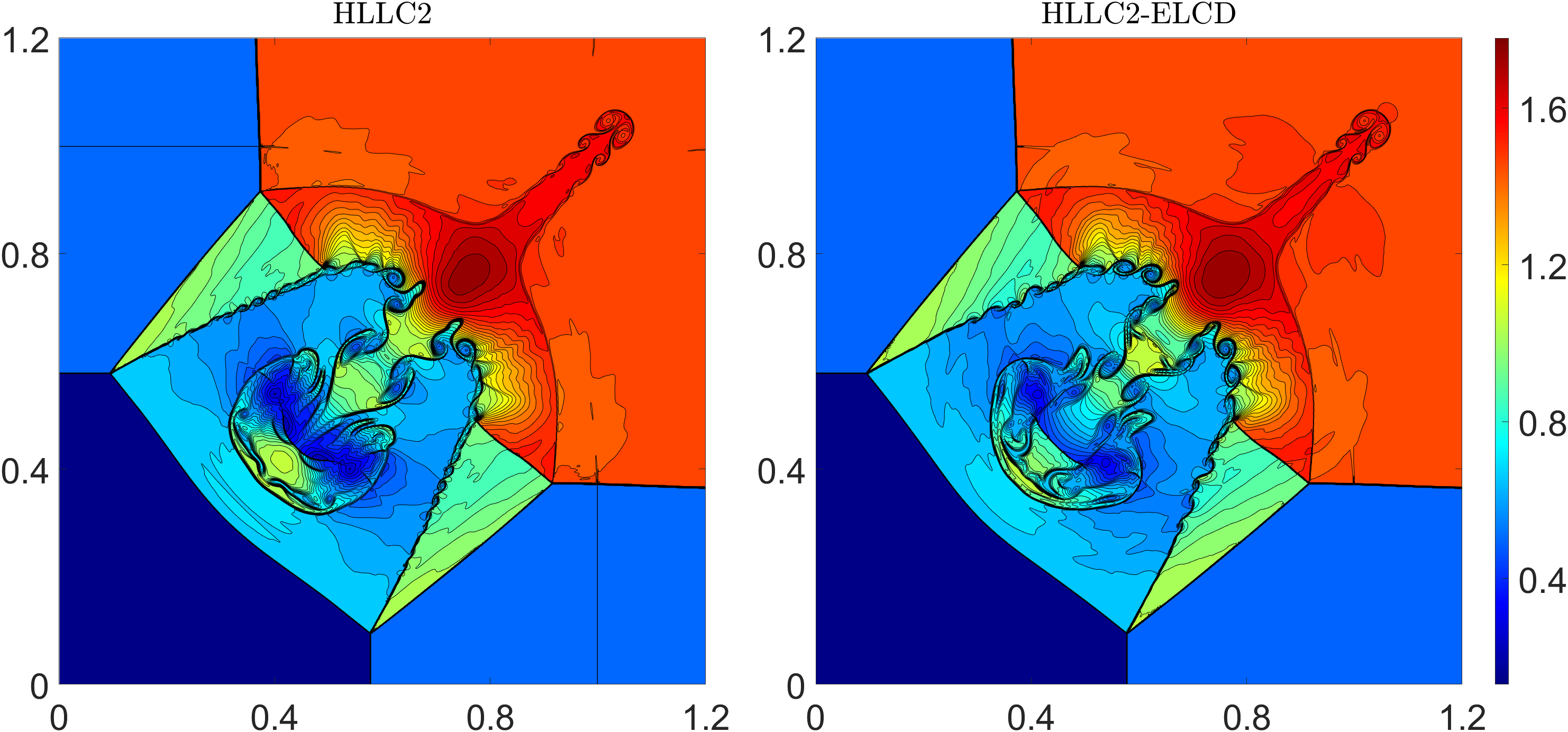}}
\vskip4pt
\centerline{\includegraphics[trim=0cm 0cm 0cm 0cm, clip, width=0.60\linewidth]{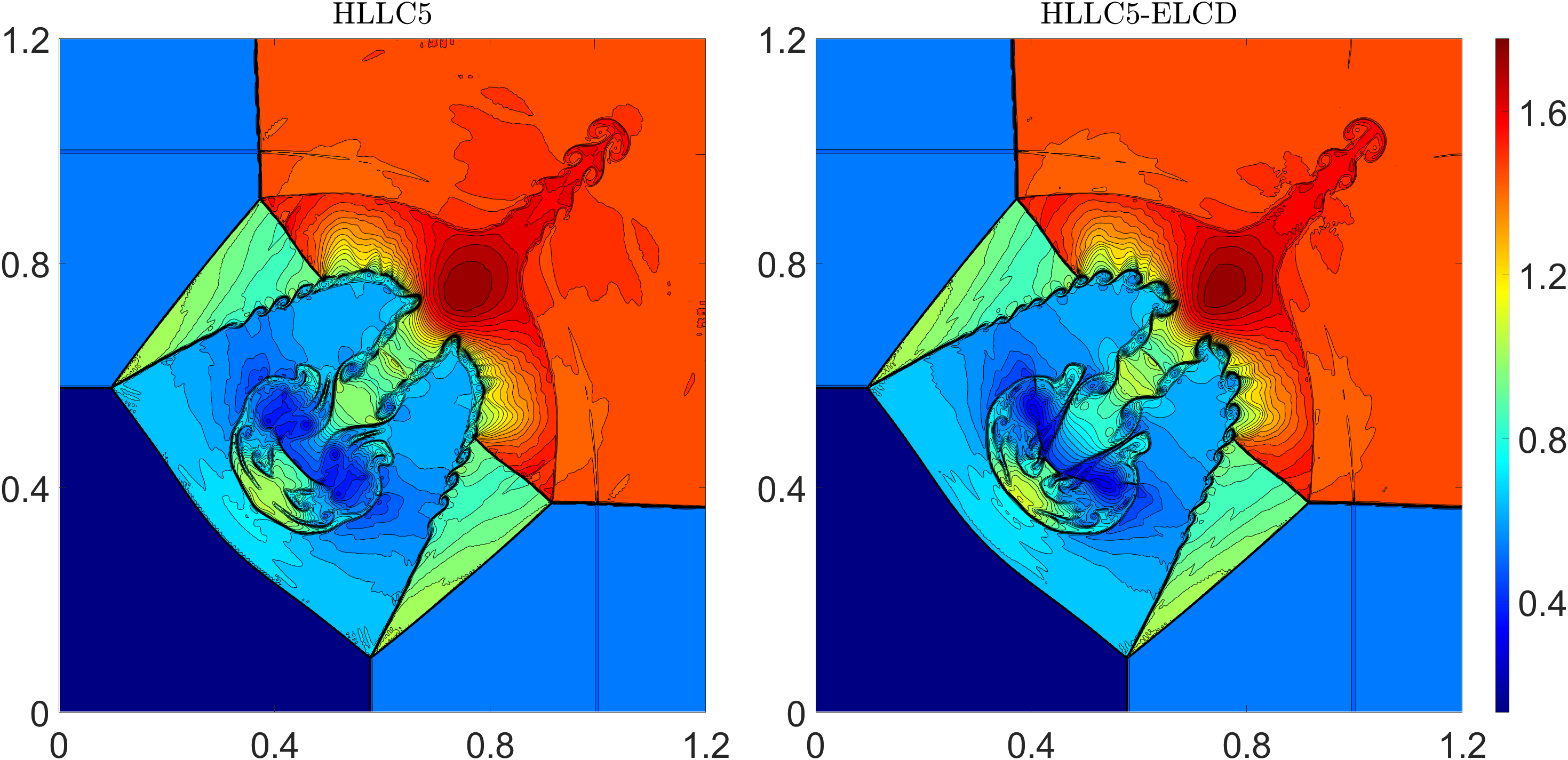}}
\caption{\sf Example 3: Density $\rho$ computed using the HLLC2 (top row) and HLLC5 (bottom row) schemes. The results obtained by the
schemes that use the ELCD procedure are plotted in the right column.\label{fig4}}
\end{figure}

We compute the numerical solution until the final time $t=1$ using the HLLC2 and HLLC5 schemes on $1000\times1000$ and $600\times600$
uniform meshes, respectively. As one can see in Fig. \ref{fig4}, the schemes, which use the ELCD procedure capture the sideband instability
of the jet more clearly than their counterparts, which use the arithmetic averages \eref{3.3a}. 
\begin{remark}
We have shown several examples, which illustrate a potential advantage of the ELCD procedure. We have tested many additional examples, both
one- and two-dimensional ones. In all of the performed numerical experiments, the use of the ELCD-based reconstructions led to either
(slightly) better results or to the results of the same quantity, which is achieved by the LCD based on the arithmetic averages of the
primitive variables.
\end{remark}

\section{Conclusion}\label{sec5}
In this paper, we have proposed an entropy-based local characteristic decomposition (ELCD), in which an average interface state is selected
by locally minimizing the entropy over several candidates constructed from the neighboring cell states. The ELCD can be readily incorporated
into existing local characteristic reconstruction/interpolation procedures. Numerical experiments with second-order FV and fifth-order FD
A-WENO schemes using four different numerical fluxes show that the resulting schemes generally resolve vortical and other complex wave
structures sharper than their counterparts, which use the arithmetic averages. We emphasize that the term ``entropy-based'' refers only to
the criterion used for selecting the average interface states and does not imply that the resulting schemes produce a smaller entropy. These
results demonstrate that the choice of the average interface states can improve the achieved resolution. 

It is important to stress that a similar ELCD linearization can be used not only for the Euler equations of gas dynamics, but for any
hyperbolic system that has convex entropy. The process of selecting the average interface states, described in Algorithm 1, should be
adjusted for the hyperbolic system at hands.


\bigskip 

\noindent{\bf Conflicts of interest.} On behalf of all authors, the corresponding author states that there is no conflict of interest.

\medskip 
\noindent{\bf Data and software availability.} The data that support the findings of this study as well as FORTRAN codes developed by the
authors and used to obtain all of the presented numerical results are available from the corresponding author upon reasonable request.

\bibliographystyle{spmpsci}      % mathematics and physical sciences
\bibliography{ref}   % name your BibTeX data base

\end{document}